# BEST MÖBIUS APPROXIMATIONS OF CONVEX AND CONCAVE MAPPINGS

MARTIN CHUAQUI AND BRAD OSGOOD

*In memory of Peter Duren*

ABSTRACT. We study the best Möbius approximations (BMA) to convex and concave conformal mappings of the disk, including the special case of mappings onto convex polygons. The crucial factor is the location of the poles of the BMAs. Finer details are possible in the case of polygons through special properties of Blaschke products and the prevertices of the mapping function.

## 1. INTRODUCTION

The *best Möbius approximation* (BMA) to a locally injective analytic function $f(z)$ at a point $\zeta$ is the unique Möbius transformation $Mf(z,\zeta)$ (in $z$) that agrees with $f(z)$ to second order at $\zeta$. Explicitly,

$$(1.1) \qquad Mf(z,\zeta) = f(\zeta) + \frac{(z-\zeta)f'(\zeta)}{1 - \frac{1}{2}(z-\zeta)\frac{f''(\zeta)}{f'(\zeta)}}.$$

This can be derived most easily by rearranging the Taylor approximation near $z = 0$ while maintaining second order contact:

$$f(z) = a_0 + a_1 z + a_2 z^2 + o(z^3) = a_0 + a_1 z \left[1 + (a_2/a_1)z\right] + o(z^3)$$

$$= a_0 + \frac{a_1 z}{1 - (a_2/a_1)z} + o(z^3).$$

A given conformal mapping, say of the unit disk $\mathbb{D}$, determines a family of BMAs as $\zeta$ varies over $\mathbb{D}$. One might ask how the properties of $f$ affect the properties of its BMAs, and vice versa. We have found it interesting to focus attention on the poles of the BMAs, exploring the connection with convex and concave functions, in particular, in the case of mappings onto the interior and exterior of convex polygons. Examples to keep in mind are the mappings onto a strip and a convex

---

The first author was partially supported by Fondecyt Grant #1190380.
*Key words*: best Möbius approximation, convex mappings, concave mappings, polygons, Blaschke product.
*2020 AMS Subject Classification*: Primary 30C45; Secondary 30C30.

sector:

$$L(z) = \frac{1}{2} \log \frac{1+z}{1-z}, \quad A(z) = \frac{1}{2\alpha}\left[\left(\frac{1+z}{1-z}\right)^\alpha - 1\right], \quad 0 < \alpha \leq 1,$$

with BMAs given by

$$ML(z,\zeta) = \frac{1}{2}\log\frac{1+\zeta}{1-\zeta} + \frac{\zeta - z}{z\zeta - 1}$$

and

$$MA(z,\zeta) = \frac{1}{2\alpha}\left(\left(\frac{1+\zeta}{1-\zeta}\right)^\alpha \frac{1 - z\zeta + \alpha(z-\zeta)}{1 - z\zeta - \alpha(z-\zeta)} - 1\right).$$

Of course, if $f$ itself is a Möbius transformation then $Mf(z,\zeta) = f(z)$. Also, for the record, we note that one can write

$$Mf(z,\zeta) = \frac{a(\zeta)z + b(\zeta)}{c(\zeta)z + d(\zeta)}$$

for

$$a(\zeta) = \frac{2(f')^2 - ff''}{2f'}, \quad b(\zeta) = f - \zeta a(\zeta), \quad c(\zeta) = -\frac{f''(\zeta)}{2f'(\zeta)}, \quad d(\zeta) = 1 + \zeta c(\zeta)$$

where $f$ and its derivatives are evaluated at $\zeta$. These quantities satisfy the normalization

$$a(\zeta)d(\zeta) - b(\zeta)c(\zeta) = f'(\zeta).$$

1.1. **Convex and concave mappings.** Recall that a function $f$ is convex in $\mathbb{D}$ if and only if

$$\operatorname{Re}\left\{1 + z\frac{f''(z)}{f'(z)}\right\} \geq 0,$$

which can be rewritten as

$$1 + z\frac{f''(z)}{f'(z)} = \frac{1 + h(z)}{1 - h(z)}$$

for some analytic $h$ with $|h| \leq 1$ in $\mathbb{D}$. Since $h(0) = 0$ we can further write

(1.2) $$1 + z\frac{f''(z)}{f'(z)} = \frac{1 + z\varphi(z)}{1 - z\varphi(z)}$$

with $\varphi$ analytic and bounded by 1 in $\mathbb{D}$. Solving for $\varphi$,

(1.3) $$\varphi(z) = \frac{f''(z)}{2f'(z) + zf''(z)}.$$

A mapping $g$ will be concave in $\mathbb{D}$ provided

$$\operatorname{Re}\left\{1 + z\frac{g''(z)}{g'(z)}\right\} \leq 0,$$

which implies that $g$ must have a (simple) pole at $z = 0$. Arguing as before, we now obtain the representation

$$1 + z\frac{g''(z)}{g'(z)} = -\frac{1 + z\omega(z)}{1 - z\omega(z)} \tag{1.4}$$

for some $\omega$ analytic and bounded by 1 in $\mathbb{D}$ that must also vanish at the origin. From here,

$$z^2\omega(z) = \frac{2g'(z) + zg''(z)}{g''(z)}. \tag{1.5}$$

It was shown in [4] that the function $\varphi$ in the representation (1.3) is a Blaschke product $B(z)$ of degree $n - 1$ precisely when $f$ maps $\mathbb{D}$ onto a convex $n$-gon. Furthermore, the prevertices are the roots of $zB(z) = 1$. Similarly, for a mapping $g$ onto the exterior of a convex $n$-gon with $g(0) = \infty$, the function $\omega$ in (1.5) is a Blaschke product $B(z)$ of degree $n - 1$ with the additional requirement that $B(0) = 0$. The prevertices correspond to the roots of the equation $zB(z) = 1$. See [1], [5] for further results on concave mappings and Blaschke products associated with arbitrary Schwarz-Christoffel mappings.

Peter Duren was a collaborator and a close friend of the authors. This paper is in line with some of our joint work and we dedicate it to his memory.

## 2. Mappings onto convex domains and their complements

We begin with the following general result on the location of the pole of the BMA for a locally injective mapping.

**Theorem 2.1.** *Let $f$ be locally injective in a domain $D$, and let $\zeta \in D$. Then the pole $p(\zeta)$ of the BMA of $f$ at $\zeta$ satisfies the following:*

*(a) The modulus of $|p(\zeta)|$ is bigger than, equal to, or smaller than $|\zeta|$ according to whether $\operatorname{Re}\{1 + \zeta(f''(\zeta)/f'(\zeta))\}$ is positive, $0$, or negative.*

*(b) The points $0, \zeta, p(\zeta)$ are colinear if and only if $\operatorname{Im}\{\zeta(f''/f')\} = 0$.*

*(c) $p(\zeta) = -\zeta$ if and only if $1 + \zeta(f''(\zeta)/f'(\zeta)) = 0$.*

*Proof.* From (1.1), the pole of $Mf(z, \zeta)$ is

$$p(\zeta) = \zeta + 2\frac{f'(\zeta)}{f''(\zeta)}, \tag{2.1}$$

which includes the case $p = \infty$ when $f''(\zeta) = 0$. Writing

$$1 + \zeta\frac{f'(\zeta)}{f''(\zeta)} = h + ik,$$

we have

$$p(\zeta) = \zeta + \frac{2\zeta}{h - 1 + ik} = -\zeta \frac{1 + h + ik}{1 - h - ik}. \tag{2.2}$$

Hence

$$|p(\zeta)|^2 - |\zeta|^2 = \frac{4h}{(1-h)^2 + k^2}, \tag{2.3}$$

which shows that $|p(\zeta)| - |z|$ has the same sign as $h$. This proves the first part.

For part (b), we see from (2.2) that $p(\zeta)$ and $\zeta$ lie on the same line through the origin iff $(1 + h + ik)/(1 - h - ik)$ is real, which occurs iff $k = 0$. We also see from this that $p(\zeta) = -\zeta$ iff $h = k = 0$, proving part (c). □

Equation (2.1) shows that $p(\zeta) = \zeta$ cannot occur for $\zeta \in D$, as it would require $f''(\zeta) = \infty$. Nevertheless, this could very well happen in the limit at some point $\zeta_0 \in \partial D$, and we see that

$$\lim_{\zeta \to \zeta_0} p(\zeta) = \zeta_0 \quad \Longleftrightarrow \quad \lim_{\zeta \to \zeta_0} \frac{f''(\zeta)}{f'(\zeta)} = \infty.$$

The figure below is a fair representation of the location of the pole relative to the base point denoted by $z$. For example, $p(z) = 0$ iff $h = -1$, $k = 0$, and $p(z) = -1$ iff $h = k = 0$. The shaded regions correspond to the different signs of $h$. The pole lies on the line through the origin and $z$ iff $k = 0$, and it will lie on the line through the origin perpendicular to this direction iff $h^2 + k^2 = 1$.

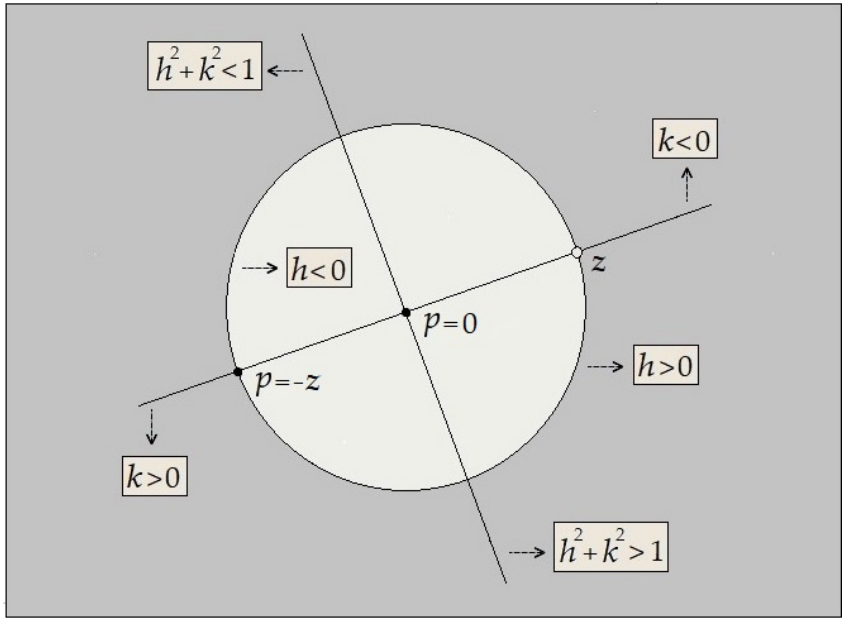

Here is Theorem 2.1 in action – more can be concluded on the location of the poles when $\operatorname{Re}\{1 + z(f''/f')\}$ does not change sign in $\mathbb{D}$.

**Theorem 2.2.** *Let $f$ be locally univalent in $\mathbb{D}$. Then $f$ is convex univalent if and only if the pole of every BMA of $f$ lies outside $\mathbb{D}$.*

The proof will show that if the pole of a single BMA lies on $\partial \mathbb{D}$ then $f$ must be a halfplane mapping, and thus, the pole of every BMA lies on $\partial \mathbb{D}$. As soon as one pole, say $\zeta_0$, lies inside the disk the image of the circle $|z| = |\zeta_0|$ will fail to be convex at the corresponding image point.

*Proof.* Let $f$ be locally injective. We see that

$$p(\zeta) = \frac{2f'(\zeta) + \zeta f''(\zeta)}{f''(\zeta)} = \frac{1}{\varphi(\zeta)}$$

for the function $\varphi$ characterizing the convexity of $f$. Therefore, $f$ is convex iff $|\varphi(\zeta)| \leq 1$ for all $|\zeta| < 1$, iff $|p(\zeta)| \geq 1$ for all $|\zeta| \leq 1$. If for a convex mapping $f$, $|p(\zeta)| = 1$ at some $\zeta \in \mathbb{D}$, then $\varphi$ must be a constant of absolute value 1, which implies that $f$ is a halfplane mapping. □

There is a corresponding result for concave mappings.

**Theorem 2.3.** *Let $g$ be locally univalent in $\mathbb{D}$ with $g(0) = \infty$. Then $g$ is concave univalent if and only if the pole of every BMA of $g$ lies in $\mathbb{D}$.*

*Proof.* The mapping $g$ is concave iff the function $\omega$ in (1.5) satisfies $|\omega| \leq 1$, which by (2.1) is equivalent to the poles of the BMAs lying inside the closed disk $\overline{\mathbb{D}}$. Since $\omega(0) = 0$ we know that $|\omega| \leq |z| < 1$, and thus all poles lie in the open disk $\mathbb{D}$. □

We introduce a correspondance between (normalized) convex and concave mappings, expressed in terms of their derivatives.

**Theorem 2.4.** *Let $f$ be locally univalent in $\mathbb{D}$ with $f''(0) = 0$, and let $g$ be defined by the conditions*

$$(2.4) \qquad f'(z)g'(z) = -\frac{1}{z^2} \ , \ g(0) = \infty \, .$$

*Then for $0 < |z| < 1$*

$$(2.5) \qquad \operatorname{Re}\left\{1 + z\frac{g''(z)}{g'(z)}\right\} = -\operatorname{Re}\left\{1 + z\frac{f''(z)}{f'(z)}\right\} \, .$$

*In particular, $f$ is convex if and only if $g$ is concave.*

*Proof.* The condition $f''(0) = 0$ is required to ensure that $g'$ will have a well-defined primitive in $\mathbb{D}\setminus\{0\}$, while the double pole of $g'$ at $z = 0$ ensures that $g$ has a simple pole there. The equation (2.5) follows at once from (2.4)                        □

We will have more to say on the condition $f''(0) = 0$ later.

Under this correspondance the total curvature along the images $f(C), g(C)$ of a circular arc $C : z = re^{it}$, $a < t < b$, $0 < r < 1$ are the same. Indeed, for the total curvature of $f(C)$ we have

$$\int_{f(C)} k(w)|dw| = \int_a^b \frac{1}{|f'|}\operatorname{Re}\left\{1 + z\frac{f''}{f'}\right\}|f'|dt = \int_a^b \operatorname{Re}\left\{1 + z\frac{f''}{f'}\right\} dt$$

$$= \int_a^b \left|\operatorname{Re}\left\{1 + z\frac{g''}{g'}\right\}\right| dt = \int_{g(C)} k(w)|dw|.$$

Thus, the angle of inclination of the tangent along both of the arcs $f(C)$ and $g(C)$ have the same total variation. The conformal mapping $g \circ f^{-1}$ of one convex domain onto the complement of another convex domain will preserve this variation.

To exemplify, we analyze briefly this correspondance for the mapping $L$ onto a parallel strip, and the mapping $K$ obtained from a Möbius shift of $A$ that ensures vanishing second derivative at the origin. The lens $K(\mathbb{D})$ consists of two circular arcs in the upper and lower halfplanes that meet at $\pm 1/\alpha$ at an angle $\alpha\pi$.

For the mapping $L$ we find that the corresponding concave mapping will have

$$g'(z) = 1 - \frac{1}{z^2},$$

hence

$$g(z) = z + \frac{1}{z}$$

maps $\mathbb{D}$ onto the complement of the interval $[-2, 2]$.

On the other hand, and despite there being an explicit expression for $K'$, it is rather complicated to find a formula for the concave counterpart. Nevertheless, from symmetry and the above discussion on total curvature, we see that $g(\partial\mathbb{D})$ must also consist of two analytics arcs which are symmetric with respect to both axes and meet at an angle of $\alpha\pi$.

3. Mappings onto convex polygons and their complements

We turn to the special case of mappings onto convex polygons and their complements.

**Theorem 3.1.** *Let $f$ map $\mathbb{D}$ onto a convex polygon and let $z_1, \ldots, z_n \in \partial\mathbb{D}$ be the prevertices, $n \geq 2$. Let $\zeta \in \partial\mathbb{D}$ be on the open arc $I_k$ from $z_k$ to $z_{k+1}$. Then the pole of the BMA to $f$ at $\zeta$ lies in the complementary open arc $\partial\mathbb{D}\setminus\overline{I_k}$. The pole* *depends on $\zeta$ in a one-to-one way and covers the entire arc $\partial\mathbb{D}\backslash\overline{I_k}$ in a clockwise sense.*

*Proof.* As mentioned earlier ([3], [4]), for a mapping $f$ onto a convex polygon, the function $\varphi$ in (1.2) is a Blaschke product $B(z)$ of degree $n-1$ and the pole of the BMA of $f$ at a point $\zeta \in I_k$ is given by

$$p(\zeta) = \frac{1}{\varphi(\zeta)} = \overline{B(\zeta)}.$$

The argument of $\overline{B(\zeta)}$ is a decreasing function of $\arg\{\zeta\}$, and since the prevertices are the roots of $zB(z) = 1$, we see that as $\zeta$ covers $I_k$, the pole $\overline{B(\zeta)}$ will cover the complementary open arc $\partial\mathbb{D}\backslash\overline{I_k}$. This dependance is one-to-one because $B' \neq 0$ on $I_k$. □

Complementing Theorem 3.1 is:

**Theorem 3.2.** *Let $g$ map $\mathbb{D}$ onto the complement of a convex polygon and let $z_1, \ldots, z_n \in \partial\mathbb{D}$ be the prevertices, $n \geq 2$. Let $\zeta \in \partial\mathbb{D}$ be on the open arc $I_k$ from $z_k$ to $z_{k+1}$. Then the pole of the BMA to $g$ at $\zeta$ covers $\partial\mathbb{D}$ moving in the positive sense, starting at $p(z_k) = z_k$ and ending at $p(z_{k+1}) = z_{k+1}$ with a total variation in argument of $2\pi + t_{k+1} - t_k > 2\pi$.*

*Proof.* The proof is almost identical to the proof of Theorem 3.1, with the difference that now

$$p(\zeta) = \zeta + 2\frac{g'(\zeta)}{g''(\zeta)} = \zeta^2 B(\zeta),$$

for a Blaschke product of degree $n-1$ (that vanishes at the origin). The theorem then follows from the fact that the prevertices are the roots of $zB(z) = 1$. □

There is an interesting observation. If $T$ is the BMA to $f$ at a point $\zeta \in \partial\mathbb{D}$ between consecutive prevertices, then $T(\mathbb{D})$ must be a disk or a halfplane containing $f(\zeta)$ on its boundary. But it must also have the same tangent line at that point, and we conclude that $f(\mathbb{D})$ is a halfplane with a boundary equal to the extended line from $f(z_k)$ to $f(z_{k+1})$. It is also the limit of the BMAs to $f$ along a sequence $\zeta_n \in \mathbb{D}$ converging to $\zeta$, so we conclude that $T(\mathbb{D})$ is the halfplane containing $f(\mathbb{D})$. The same argument applies to any convex or concave mapping $f$ that is $C^2$ up to the boundary: if the curvature of $f(\partial\mathbb{D})$ vanishes at some boundary point $f(\zeta)$ then the BMA to $f$ at $\zeta$ maps $\mathbb{D}$ to a supporting halfplane containing $f(\mathbb{D})$.

We finish with the following theorems on the behavior of $|f'|$ and $|g'|$ on the open arcs $I_k$. Let $z_k = e^{it_k}$, for $0 \leq t_1 < t_2 < \ldots < t_n < 2\pi$.

**Theorem 3.3.** *Let $f$ map $\mathbb{D}$ onto a convex $n$-gon with open arcs $I_k \subset \partial\mathbb{D}$ determined by the prevertices. Then $u(t) = |f'(e^{it})|^{-1/2}$ is a concave function on each open arc $I_k$. In particular, $|f'(e^{it}|$ is a convex function on the $I_k$.*

*Proof.* Consider the curve $\phi(t) = f(e^{it})$ as a function of $t = \arg\{z\}$, $z \in I_k$. Then $u = |\phi'(t)|^{-1/2} = |f'(e^{it})|^{-1/2}$ is a positive solution of

$$u'' + \frac{1}{2}\mathrm{Re}\{S\phi\}u = 0,$$

where $S\phi$ is the Schwarzian derivative of $\phi$. From the chain rule for the Schwarzian we find that $S\phi = (ie^{it})^2 Sf + S\sigma$, where $\sigma(t) = e^{it}$. It follows that

$$S\phi = -z^2 Sf + \frac{1}{2}, \quad z = e^{it},$$

and from (1.2) with $\varphi = B$ a Blaschke product of degree $n - 1$, we find that

$$z^2 Sf = \frac{2z^2 B'}{(1 - zB)^2} = \frac{2|B'|w}{(1-w)^2}, \quad w = zB.$$

Here we have used that $zB'/B = |B'| > 0$ on $\partial \mathbb{D}$. Therefore, $z^2 Sf$ has negative real part on $I_k$ because $w/(1-w)^2$ is the Koebe mapping. We conclude that $\mathrm{Re}\{S\phi\} > 0$, hence $u$ is concave. From this, $|f'| = u^{-2}$ is convex. □

Again there is a complementary result.

**Theorem 3.4.** *Let $g$ map $\mathbb{D}$ onto the complement of a convex $n$-gon with $g(0) = \infty$ and open arcs $I_k \subset \partial \mathbb{D}$ determined by the prevertices. Then $v(t) = |g'(e^{it})|^{-1/2}$ is a convex function on each open arc $I_k$.*

*Proof.* Following the same argument as above with $\psi(t) = g(e^{it})$ and $v = |g'|^{-1/2}$, we will show that $-z^2 Sg + \frac{1}{2}$ has negative real part on each $I_k$. From (1.5) we find that

$$z^2 Sg = -\frac{2z^2 B'}{(1 - zB)^2} - \frac{4zB}{(1 - zB)^2} = -\frac{2|B'|w}{(1-w)^2} - \frac{4w}{(1-w)^2},$$

hence

$$S\psi = \frac{2|B'|w}{(1-w)^2} + \frac{4w}{(1-w)^2}.$$

The first term on the right hand-side has negative real part, while the real part of the second term is smaller than $-1$. This shows that $\mathrm{Re}\{S\psi\} < 0$, and thus $v$ is convex as claimed. □

It follows from Theorems 3.3, 3.4 that on each arc $I_k$ there exists a unique point $w_k$ where $|f'|$ and $|g'|$ attain, respectively, their maximum and minimum value in $I_k$. These points admit the following interpretations. First, it follows from Theorem 2.1, part (c), that such points $w_k$ correspond to when $p(\zeta) = -\zeta$. Secondly, it can be shown that the correspondence under Alexander's theorem of a mapping onto a convex polygon is a starlike mapping onto the plane minus the same number of slits to infinity, and furthermore, that for the Blaschke product arising from the convex mapping, the roots of $zB(z) = -1$ correspond to the pre-images of the starlike mapping of the finite ends of the slits ([2]). The correspondance given in Theorem 2.4 will associate a mapping $f$ onto a convex $n$-gon with a concave mapping ommiiting a certain other $n$-gon, with the property that prevertices and exterior angles will be the same. For the special case when $n = 3$, the two triangles must therefore be similar. Hence, we draw the curious fact that there must be a constant $C > 0$ such that for $k = 1, 2, 3$

$$\int_{I_k} |f'| dt = C \int_{I_k} |g'| dt = C \int_{I_k} \frac{1}{|f'|} dt\,.$$

**Lemma 3.5.** *Let $\tau : \partial \mathbb{D} \to \mathbb{R}$ be a continuous function that is strictly increasing as a function of the $x$-coordinate. Let $z_1, z_2, z_3 \in \partial \mathbb{D}$ be three distinct points ordered according to increasing argument, and denote by $J_k \subset \partial \mathbb{D}$ the closed arc joining $z_k$ to $z_{k+1}$, $k = 1, 2, 3$, with the understanding that $z_4 = z_1$. Then for some $i \neq k$,*

$$\max_{J_i} \tau \leq \min_{J_k} \tau\,.$$

The lemma is basically claiming that for some $x_0 \in (-1, 1)$, there are two arcs $J_i, J_k$ that lie on different sides of the line $x = x_0$.

*Proof.* After relabelling if necessary, either $\operatorname{Re}\{z_1\} \leq \operatorname{Re}\{z_2\} \leq \operatorname{Re}\{z_3\}$ or $\operatorname{Re}\{z_1\} \leq \operatorname{Re}\{z_3\} \leq \operatorname{Re}\{z_2\}$. Both cases are treated similarly, and we will assume the former. Notice also that there can be only one equality between these quantities. Assume there is one such equality, say $\operatorname{Re}\{z_1\} = \operatorname{Re}\{z_2\}$, meaning that $z_1, z_2$ are conjugate and that $\operatorname{Re}\{z_1\} < \operatorname{Re}\{z_3\}$. In this case, $\max_{J_1} \tau \leq \min_{J_2} \tau$. If, on the other hand, $\operatorname{Re}\{z_1\} < \operatorname{Re}\{z_2\} < \operatorname{Re}\{z_3\}$, then once again $\max_{J_1} \tau \leq \min_{J_2} \tau$. This finishes the proof. □

In the following theorem, the condition $f''(0) = 0$ implies that $f(0)$ is a critical point of the Poincaré density of the triangle. Since $f$ is convex, and thus belongs the Nehari class ([10], [9], [8], [3]), a critical point can only occur if the triangle is bounded ([7]). From the classical expression for $f''/f'$ as a rational function, this normalization corresponds to the equation

$$\alpha_1 z_1 + \alpha_2 z_2 + \alpha_3 z_3 = 0$$

involving the prevertices and the exterior angles, and which determines the prevertices up to rotation. We finish with the following characterization.

**Theorem 3.6.** *Let $f$ map $\mathbb{D}$ onto a bounded triangle, with prevertices $z_1, z_2, z_3$ and subarcs on $\partial \mathbb{D}$ joining them, $I_1, I_2, I_3$. Then the following conditions are equivalent:*

*(a) $f''(0) = 0$;*

*(b) there exists a positive constant $C$ such that for $k = 1, 2, 3$*

(3.1) $$\int_{I_k} |f'| dt = C \int_{I_k} \frac{1}{|f'|} dt\,.$$

*Proof.* The discussion preceeding the lemma shows that $(a)$ implies $(b)$. For the reverse implication, consider a mapping $f$ onto the triangle with $f''(0) = 0$. Since any other conformal mapping onto the triangle is of the form $g = f \circ \sigma$ for some autmorphism $\sigma$ of the disk, it suffices to show that $g$ cannot satisfy (0.1) unless $\sigma(0) = 0$. Suppose that $\sigma(0) \neq 0$, say without loss of generality, that $\sigma(0) = -a \in (-1, 0)$. From the expression $\sigma(z) = c(z+a)/(1+az)$, $|c| = 1$, it follows that

$$|\sigma'(e^{it})| = \frac{1-a^2}{1+a^2+2a\cos(t)},$$

hence $|\sigma'|$ satisfies the hypothesis in Lemma 0.1. Let $w_k = \sigma^{-1}(z_k)$ and $J_k = \sigma^{-1}(I_k)$. The change of variables $e^{it} = \sigma(e^{is})$ gives

$$\int_{I_k} |f'| dt = \int_{J_k} |g'| ds,$$

while the same susbstitution and the mean value theorem give

$$\int_{I_k} \frac{1}{|f'|} dt = \int_{J_k} \frac{|\sigma'|^2}{|g'|} ds = c_k \int_{J_k} \frac{1}{|g'|} ds,$$

where $c_k$ lies between $\min_{J_k} |\sigma'|^2$ and $\max_{J_k} |\sigma'|^2$. Thus,

$$\int_{J_k} \frac{1}{|g'|} ds = \frac{1}{c_k} \int_{I_k} \frac{1}{|f'|} dt = \frac{1}{Cc_k} \int_{I_k} |f'| dt = \frac{1}{Cc_k} \int_{J_k} |g'| ds,$$

which together with Lemma 3.5 shows that $g$ cannot satisfy condition (b) of the theorem. This finishes the proof. $\square$

We finish with an observation that is derived from the proof of Theorem 3.6. If $f, g$ are conformal mappings onto the same bounded triangle with $f''(0) = 0$, then there exists a constant $C > 0$ such that for $k = 1, 2, 3$

$$Cb^{-2} \leq \int_{J_k} \frac{1}{|g'|} ds \leq \int_{J_k} |g'| ds \leq Cb^2 \int_{J_k} \frac{1}{|g'|} ds,$$

where $b = (1+a)/(1-a)$ and $a = (f^{-1} \circ g)(0)$. This follows from the fact that $b, b^{-1}$ are the maximum and minimum values of $|\sigma'|$ on $\partial \mathbb{D}$, as considered in the proof.

## References

[1] B. Bhowmik, S. Ponnusamy, and K.-J. Wirths, *Concave functions, Blaschke products, and polygonal mappings*, Siberian Math. J. 50 (2009), 609–615. MR2583615 (2011a:30123)

[2] M. Chuaqui and B. Osgood, *On convolution, convex and starlike mappings*, Stud. Babeş-Bolyai Math. 67 (2022), 431-440.

[3] M. Chuaqui, P. Duren and B. Osgood, *Schwarzian derivatives of convex mappings*, Ann. Acad. Scie. Fenn. Math. 36 (2011), 1-12

[4] M. Chuaqui, P. Duren and B. Osgood, *Concave conformal mappings and pre-vertices of Schwarz-Christoffel mappings*, Proc. AMS 140 (2012), 3495-3505

[5] M. Chuaqui and Ch. Pommerenke, *On Schwarz-Christoffel mappings*, Pacific J. Math. 270 (2014), 319-334

[6] P. Duren, *Univalent Functions*, Grundlehreen der mathematischen Wissenschaften 259, Springer-Verlag, New York, 1983

[7] F.W. Gehring and Ch. Pommerenke, *On the Nehari univalence criterion and quasicircles*, Comm. Math. Helv. 59 (1984), 226-242.

[8] S.-A. Kim and D. Minda, *The hyperbolic and quasihyperbolic metrics in convex regions*, J. Analysis 1 (1993), 109–118.

[9] W. Koepf, *Convex functions and the Nehari univalence criterion*, in Complex Analysis, Joensuu 1987 (I. Laine, S. Rickman, and T. Sorvali, editors), Lecture Notes in Math., vol. 1351, Springer-Verlag, Berlin, 1988, pp. 214–218.

[10] Z. Nehari, *A property of convex conformal maps*, J. Analyse Math. 30 (1976), 390-393

Facultad de Matemáticas, Pontificia Universidad Católica de Chile, Santiago, Chile, mchuaqui@mat.uc.cl

Electrical Engineering, Stanford University, Palo Alto, CA 94025, USA, osgood@ee.stanford.edu